\documentclass{amsart}

\usepackage{xcolor}
\definecolor{pagecolor}{gray}{1}

\usepackage[pdfborder={0 0 0}, draft=false, pagebackref, unicode=true]{hyperref}
\usepackage{enumitem, esint, amssymb, mathrsfs, mathtools, amsthm, tikz, enumitem}
\usetikzlibrary{shadings}

\newtheorem{thm}[equation]{Theorem}
\theoremstyle{definition}
\newtheorem{defn}[equation]{Definition}
\newtheorem{rmk}[equation]{Remark}

\newcommand\corkscrew{\varkappa}

\newcommand\dist{\mathop{\mathrm{dist}}\nolimits}
\newcommand\Div{\mathop{\mathrm{div}}\nolimits}
\newcommand\diam{\mathop{\mathrm{diam}}\nolimits}
\newcommand\R{\mathbb{R}}

\newcommand\dmn{{\mathfrak{n}}}  
\newcommand\pdmn{{\mathfrak{n}}} 
\newcommand\bdmn{{\mathfrak{d}}} 
\newcommand\dmnMinusOne{{\mathfrak{n}-1}}
\newcommand\pdmnMinusOne{{(\mathfrak{n}-1)}}
\newcommand\pbdmn{{\mathfrak{d}}}
\newcommand\bpdmn{{\mathfrak{d}}}

\makeatletter\def\HyPsd@CatcodeWarning#1{}\makeatother

\title{The Poisson-Dirichlet problem in domains with Ahlfors regular boundary}

\author{Ariel Barton, Svitlana Mayboroda, and Alberto Pacati}

\begin{document}

\subjclass[2020]{35J25 (Primary) 35A01, 35A02 (Secondary)}

\begin{abstract} We present an announcement of some recent results concerning well-posedness of the Poisson-Dirichlet problem with boundary data in Besov spaces with fractional smoothness. This is a far-reaching generalization as previously known theorems concerning well-posedness of the Poisson problem in such intermediate smoothness classes were mostly restricted to the context of Lipschitz domains and coefficients satisfying strong regularity assumptions.
\end{abstract}

\maketitle

\markright{THE POISSON-DIRICHLET PROBLEM IN DOMAINS WITH AHLFORS BOUNDARY}

\section{Introduction: the general context}

The purpose of this note is to describe some recent results concerning the Poisson problem, and the Dirichlet problem with data of fractional smoothness, for second order elliptic operators with real coefficients in very general domains.

Poisson's equation $-\Delta u=h$, for a given function~$h$ and an unknown function~$u$, was historically one of the first differential equations to be studied, due in large part to its many applications in physics and beyond. A natural generalization is to replace the Laplace operator $-\Delta$ by an arbitrary divergence form elliptic operator $-\Div A\nabla$ for some specified matrix $A$ (which may be only bounded, measurable, and uniformly positive definite). Considering Poisson's equation not in all of Euclidean space, but only in a given open subset~$\Omega$, requires specifying the behavior of~$u$ on the boundary~$\partial\Omega$; we will impose Dirichlet boundary conditions, that is, specify the values of $u$ on the boundary. (It is also common in the theory to specify Neumann boundary conditions, that is, to specify $\nu\cdot A\nabla u$, where $\nu$ denotes the unit outward normal vector to~$\Omega$; we refer the interested reader to \cite{FabMM98,Zan00,BarM16A,AmeA18} for some results on the Poisson-Neumann problem.)

We thus will discuss the Poisson-Dirichlet problem
\begin{equation}\label{eqn:Dirichlet:Poisson:introduction} -\Div A\nabla u=h\text{ in }\Omega,\quad
u=f\text{ on }\partial\Omega
\end{equation}
for specified $\Omega$, $A$, $f$, and~$h$.

In a few situations, the range of spaces allowing for well-posedness of the Poisson-Dirichlet boundary problem above is very well understood. Such situations include $VMO$ and small-$BMO$ coefficients and domains (see \cite{MazMS10,YanYY22}, among others), and Lipschitz domains with coefficients satisfying a Dini-type continuity condition \cite{MitT06} or independent of the direction transverse to the boundary \cite{BarM16A}.

In recent years, end point results, pertaining to boundary data $f$ lying in $L^p(\partial\Omega)$, have been established considerably beyond Lipschitz domains: first, in chord-arc domains
\cite{DavJ90,JerK82}, then in so-called 1-sided non-tangentially accessible domains with a uniformly rectifiable boundary
\cite{HofM14}, and finally, on domains with a uniformly rectifiable boundary that satisfy a weak local John condition \cite{AzzHMMT20}.
The latter is a sharp condition, necessary and sufficient for well-posedness results for the Laplacian \cite{AzzHMMT20}.
Finally, boundary value problems in domains with lower dimensional boundaries were treated in \cite{DavM22} (a later but much simpler treatment is in \cite{Fen22}).

Most of these results were proven in (or may be generalized to) the generality of coefficients satisfying variations of the so-called DKP condition, which controls oscillations of the coefficients near the boundary.

Note that $t$-independence is meaningless in domains that lack a distinguised direction (in particular, in domains where there is no constant direction transverse to the boundary), and that all of these “end point” estimates are in spaces with integer smoothness. Our results, by contrast with the end point estimates listed above, need neither $t$-independence nor the DKP condition.

The results of this note (described in Sections~\ref{sec:Poisson} and~\ref{sec:extrapolation}) treat well-posedness of boundary problems in intermediate spaces of fractional smoothness. We will show, for the first time, that well posedness in a wide range of Besov spaces
can be established by a ``blind" mechanism, in the full generality of domains with Ahlfors regular boundary (possibly of lower codimension) and satisfying the weak local John condition, and general real elliptic equations, without any additional restrictions on the coefficients.
Furthermore, given end point (that is, $L^p(\partial\Omega)$) results such as those above, the range of Besov spaces can be expanded in, again, a “blind” fashion, with no further requirements on the domains or coefficents beyond those above and those required by the end point results.

The methods are new, as \cite{MayMit04A, BarM16A} relied on layer potentials, not available in such a general setting, \cite{JerK95}
worked with the Laplacian on Lipschitz domains and the case $p>1$ only, and \cite{MitT06,MazMS10,YanYY22} were confined to continuous, $VMO$, or near-$VMO$ coefficients and Lipschitz domains. The full proofs will appear in \cite{BarMP25pA} (the case $f=0$, that is, $u=0$ on~$\partial\Omega$) and \cite{BarMP25pB} (full generality) and in this note we only announce the upcoming results and put them into historic context.

We remark that $L^p(\partial\Omega)$ represents the ``zero smoothness'' end point; an alternative end point is the ``smoothness one'' end point, often called the regularity problem. The regularity problem, like the $L^p$-Dirichlet problem, has been investigated extensively. However, our results will not make use of the regularity problem, only the $L^p$-Dirichlet problem and a duality argument, to expand the range of Besov spaces. For this reason we will not discuss known results for the regularity problem.

\section{Definitions}\label{sec:dfn}

In this section we will provide all of the definitions needed to state our main results. An experienced reader can skip them and move on to Section~\ref{sec:history}. Some additional definitions are necessary to state the corollaries to our second main result, and also to discuss the history of the Poisson problem and the Dirichlet problem with boundary data in fractional smoothness spaces; these definitions will be deferred to Section~\ref{sec:dfn:2} or to the discussion.

\subsection{Geometric conditions}\label{sec:dfn:geometry}

\begin{defn}[Ahlfors regularity]\label{dfn:adr}If $\Omega\subset\R^\dmn$ and $0<\bdmn<\dmn$, then we say that the boundary $\partial\Omega$ is $\bdmn$-Ahlfors regular if there exists a constant $\mathcal{A}>0$ such that
\begin{equation*}\frac{1}{\mathcal{A}}r^\bdmn \leq \sigma(B(x,r)\cap\partial\Omega) \leq \mathcal{A}r^\bdmn\end{equation*}
for all $x\in\partial\Omega$ and all $r$ with $0<r<\diam(\Omega)$; here $\sigma$ denotes the Hausdorff $\bdmn$-dimensional surface measure.
\end{defn}

\begin{rmk}
Following \cite{DavFM21}, we do not require the boundary dimension $\bdmn$ to be an integer. Our main results (Theorems~\ref{thm:Dirichlet} and~\ref{thm:Dirichlet:Lp}) are known only in the case $\bdmn\leq\dmnMinusOne$, that is, the classical codimension~$1$ case and the higher codimension case. The lower codimensional case $\bdmn>\dmnMinusOne$ remains open.
\end{rmk}

\begin{defn}[Interior corkscrew condition]\label{dfn:iCS}
We say that $\Omega$ satisfies the {interior corkscrew condition} if there is a constant $\corkscrew\in (0,1)$ such that, if $x\in \partial\Omega$ and $0<r<\diam\Omega$, then there is a point $A_r(x)$ (called a corkscrew point) such that
\begin{equation*}
A_r(x)\in B(x,r)\quad\text{and}\quad B(A_r(x), \corkscrew r)\subset\Omega.
\end{equation*}
\end{defn}

\begin{defn}[Weak local John domain]\label{dfn:local:John}
Let $\Omega\subset\R^\dmn$ be open.
Let $x\in\partial\Omega$ and let $y\in\Omega$. We say that $\zeta$ is a {$\lambda$-carrot path} from $x$ to $y$ if $\zeta$ is a connected rectifiable path, which we require to be parameterized by arc length, such that $\zeta:[0,T]\to\Omega\cup\{x\}$ is continuous, $\zeta(0)=x$, $\zeta(T)=y$, and $\dist(\zeta(t),\partial\Omega)\geq \lambda t$ for all $t\in [0,T]$.

Now, suppose in addition that $\partial\Omega$ is $\bdmn$-Ahlfors regular for some $0<\bdmn<\dmn$.
Let $y\in\Omega$, let $\theta$, $\lambda\in (0,1]$, and let $N\geq 2$. We say that $y$ is a $(\theta,\lambda,N)$-weak local John point if there is a Borel set $F\subset B(y,N\dist(y,\partial\Omega))\cap\partial\Omega$ such that $\sigma(F)\geq \theta\sigma( B(y,N\dist(y,\partial\Omega))\cap\partial\Omega )$ and such that, if $x\in F$, then there is a $\lambda$-carrot path from $x$ to~$y$.

We say that $\Omega$ is a {weak local John domain} if there exist constants $\theta\in (0,1]$, $\lambda\in (0,1]$, and $N\geq 2$ such that every $y\in\Omega$ is a $(\theta,\lambda,N)$-weak local John point.
\end{defn}

\begin{rmk}\label{rmk:higher:John:corkscrew}A domain with lower-dimensional Ahlfors regular boundary (that is, with $\bdmn<\dmnMinusOne$) necessarily satisfies the interior corkscrew condition (see \cite[Lem\-ma~11.6]{DavFM21}). It also necessarily satisfies the weak local John condition: by \cite[Lem\-ma~2.1]{DavFM21}, such a domain must satisfy the Harnack chain condition (defined in many places in the literature, or see Definition~\ref{dfn:Harnack} below), and by connecting corkscrew points with Harnack chains we can construct carrot paths to a point $y$ in the domain from \emph{any} nearby boundary point~$x$.
\end{rmk}

\subsection{Functions and function spaces}\label{sec:dfn:function}

\begin{defn}
\label{dfn:Slobodeckij}
Let $\Gamma\subseteq\R^\dmn$ be a $\pbdmn$-Ahlfors regular set, let $0<s<1$, and let $1\leq p\leq \infty$. The Slobodeckij space $\dot \Lambda^{s,p}(\Gamma)$ is defined to be the set of all equivalence classes (up to addition of constants and sets of measure zero) of locally integrable functions $f:\Gamma\to\R$ such that the norm
\begin{equation*}\|f\|_{\dot \Lambda^{s,p}(\Gamma)}
=
\biggl(\int_{\Gamma}\int_\Gamma
\frac{|f(x)-f(y)|^p}{|x-y|^{\bdmn+ps}}
\,d\sigma(y)\,d\sigma(x)\biggr)^{1/p}
\end{equation*}
(or the corresponding $\sup$ norm when $p=\infty$) is finite.
\end{defn}

It is clear that the Slobodeckij spaces are Banach spaces.

\begin{defn}[{\cite[Definition~2.1]{HanLY99}; \cite[Definition~7]{HanY03}}]
\label{dfn:atom}
Let $\Gamma \subset\R^\dmn$ be a $\pbdmn$-Ahlfors regular set, let $z\in\Gamma$, and let $\varrho>0$.
We let $\Delta(z,\varrho)=\Gamma\cap B(z,\varrho)$.
A measurable function $a:\Gamma \to\R$ is said to be a
\emph{block in $\Delta(z,\varrho)$} if
\begin{itemize}
\item $a$ is supported in $\Delta(z,\varrho)$,
\item $|a(x)|\leq \varrho^{-\pbdmn/2}$ for all $x\in \Delta(z,\varrho)$ (and thus for all $x\in\Gamma$),
\item $|a(x)-a(y)|\leq |x-y|\varrho^{-\pbdmn/2-1}$ for all $x$, $y\in\Gamma$.
\end{itemize}
\end{defn}

\begin{rmk}In the setting of Hardy spaces, it is customary to work with atoms, which are blocks with the additional requirement $\int_{\Gamma } a\,d\sigma=0$. However, the Besov spaces are closer in character to the Hardy-Sobolev spaces, that is, functions whose gradient lies in a Hardy space. The moment condition $\int_{\R^\bdmn} \nabla a=0$ is trivially true for any compactly supported function~$a$, and so moment conditions are not usually required of Hardy-Sobolev atoms; thus, we do not require them of our blocks.
\end{rmk}

\begin{defn}
\label{dfn:atom:space}
Let $\Gamma \subset\R^\dmn$ be a $\pbdmn $-Ahlfors regular set,
let $0<s\leq 1$, and let $\pbdmn/(\bdmn+s)<p\leq 1$ (so $s+\bdmn-\bdmn/p>0$). We define the atomic space $\dot A^{s,p}(\Gamma)$ as
\begin{multline*}\dot A^{s,p}(\Gamma)=\Bigl\{\sum_{k=1}^\infty \lambda_k a_k: \sum_{k=1}^\infty |\lambda_k|^p \varrho_k^{\bdmn-\pbdmn p/2-sp}<\infty\\
\text{where each $a_k$ is a block in $\Delta(z_k,\varrho_k)$ for some $z_k\in\Gamma$}\Bigr\}
\end{multline*}
and let the norm in $\dot A^{s,p}(\Gamma)$ be given by
\begin{multline*}\|f\|_{\dot A^{s,p}(\Gamma)}=\inf\Bigl\{\Bigl(\sum_{k=1}^\infty |\lambda_k|^p \varrho_k^{\bdmn -\pbdmn p/2-sp}\Bigr)^{1/p}:\\
f=\sum_{k=1}^\infty \lambda_k a_k\text{ for
blocks $a_k$ in $\Delta(z_k,\varrho_k)$}\Bigr\}.\end{multline*}
\end{defn}
It is straightforward to establish (and we will do so in the forthcoming paper \cite{BarMP25pB}) that if $\sum_{k=1}^\infty |\lambda_k|^p \varrho_k^{\bdmn -\pbdmn p/2-sp}<\infty$, then $\sum_{k=1}^\infty \lambda_k a_k$ converges in the Banach space $\dot \Lambda^{s+\bdmn-\bpdmn/p,1}(\Gamma)$, and in fact $\dot A^{s,p}(\Gamma)$ embeds in $\dot \Lambda^{s+\bdmn-\bpdmn/p,1}(\Gamma)$. Thus the infinite sums above are meaningful.

\begin{defn}\label{dfn:Besov}
Let $\Gamma \subset\R^\dmn$ be a $\pbdmn $-Ahlfors regular set for some $0<\bdmn<\dmn$, let $0<s\leq 1$, and let $\pbdmn/(\bdmn+s)<p\leq \infty$.

We define the Besov space $\dot B^{p,p}_{s}(\Gamma)$ by
\begin{equation*}\dot B^{p,p}_{s}(\Gamma)=\begin{cases}
\dot \Lambda^{s,p}(\Gamma), & 1\leq p\leq \infty,\\
\dot A^{s,p}(\Gamma), & \frac{\pbdmn}{\bdmn+s}<p<1\end{cases}\end{equation*}
with the norm in $\dot B^{p,p}_{s}(\Gamma)$ given by the norm in the Slobodeckij space $\dot \Lambda^{s,p}(\Gamma)$ or atomic space~$\dot A^{s,p}(\Gamma)$.
\end{defn}

In Euclidean spaces, the Besov spaces $\dot B^{p,q}_{s}(\R^\bdmn)$ and the related Triebel-Lizorkin spaces $\dot F^{p,q}_{s}(\R^\bdmn)$ (they are equal if $p=q$) encompass many known function spaces. See \cite[Sections~2.3.5, 5.23, 5.2.4]{Tri83} and \cite[Proposition~3]{MenM00}. For example, the Lebesgue spaces satisfy $L^p(\R^\bdmn)=\dot F^{p,2}_{0}(\R^\bdmn)$ if $1<p<\infty$; the Sobolev spaces $\dot W^{k,p}(\R^\bdmn)$ (of functions $f$ whose $k$th gradient $\nabla^k f$ lies in $L^p(\R^\bdmn)$) satisfy $\dot W^{k,p}(\R^\bdmn)=\dot F^{p,2}_{k}(\R^\bdmn)$; the Hardy spaces satisfy $H^p(\R^\bdmn)=\dot F^{p,2}_{0}(\R^\bdmn)$ if $\frac{\bdmn}{\bdmn+1}<p<\infty$; the spaces $\dot C^\alpha(\R^\bdmn)$ of H\"older continuous functions satisfy $\dot C^\alpha(\R^\bdmn)=\dot B^{\infty,\infty}_{\alpha}(\R^\bdmn)$ for $0<\alpha<1$; the Slobodeckij spaces $\dot\Lambda^{s,p}$ satisfy the above relation $\dot B^{p,p}_{s}(\R^\bdmn)=\dot\Lambda^{s,p}(\R^\bdmn)$ for $0<s<1$ and $1\leq p\leq\infty$ (and indeed there exists a Slobodeckij type characterization of $\dot B^{p,q}_{s}(\R^\bdmn)$ for any $0<s<1$, any $1\leq p<\infty$, and $1\leq q\leq\infty$, not only for $q=p$), and the Bessel potential spaces $L^p_s=(I-\Delta)^{-s/2}(L^p)$ satisfy $L^p_s=F^{p,2}_s$. Even with the restrictions $p=q$, $0<s<1$, and $\frac{\bdmn}{\bdmn+s}<p\leq \infty$, there are many equivalent characterizations of the Besov spaces $\dot B^{p,p}_{s}(\R^\bdmn)$.

In \cite{HanS94}, Besov spaces were generalized from the Euclidean setting to the setting of spaces of homogeneous type. An active area of research consists of studying these spaces; we mention in particular that by \cite{MulY09}, our Slobodeckij characterization of the Besov spaces for $p\geq 1$ coincides with the characterization in \cite{HanS94,HanMY08}, and by \cite{HanY03}, the \emph{inhomogeneous} Besov spaces have an atomic characterization; it is natural to expect the homogeneous Besov spaces to have an atomic characterization (specifically, the atomic characterization used above), but we have not been able to find this result in the literature on spaces of homogeneous type.

\subsection{Boundary values}

\begin{defn}[Boundary trace]
Following \cite{JonW84}, if $\Omega\subset\R^\dmn$, if $\partial\Omega$ is $\bdmn$-Ahlfors regular for some $0<\bdmn<\dmn$, if $u:\Omega\to \R$ is locally integrable, and if $f:\partial\Omega\to \R$, then we say that $u=f$ on $\partial\Omega$ if
\begin{equation}\label{eqn:Tr}\lim_{r\to 0^+} \frac{1}{r^\dmn} \int_{B(x,r)\cap\Omega} |u-f(x)|=0\end{equation}
for $\sigma$-almost every $x\in\partial\Omega$.
\end{defn}
In most reasonably behaved domains, this coincides with most other notions of trace. We will use this definition in Theorems~\ref{thm:Dirichlet} and~\ref{thm:Dirichlet:Lp}. We mention that some of the results we will describe in Section~\ref{sec:history} were proven under different notions of trace (most prominently traces in the sense of Sobolev spaces); however, the different notions of trace can generally be shown to be equivalent as needed.

We observe that if $\Omega$ satisfies the interior corkscrew condition, then the value of $f$ is uniquely determined in that, for each $x\in\partial\Omega$, there can be at most one $f(x)\in\R$ such that formula~\eqref{eqn:Tr} is valid.

\subsection{Elliptic differential operators}
\begin{defn}[Ellipticity] Let $\Omega\subset\R^\dmn$ and suppose that $\partial\Omega$ is $\bdmn$-Ahlfors regular. Let $A$ be a matrix-valued function.
We say that the operator $-\Div A\nabla$ is elliptic, or, respectively, the matrix $A$ is elliptic, if each $A(x)$ is a $\dmn\times\dmn$ matrix of real, not necessarily symmetric coefficients that satisfy the following positive definiteness and boundedness conditions:
\begin{align}
\label{eqn:elliptic:introduction}
{\xi}\cdot A(x)\xi
&\geq
\lambda\dist(x,\partial\Omega)^{\bpdmn+1-\pdmn}\|\xi\|^2
,\\
\label{eqn:elliptic:bounded:introduction}
|\eta\cdot A(x)\xi|
&\leq \Lambda\dist(x,\partial\Omega)^{\bpdmn+1-\pdmn}\|\xi\|\,\|\eta\| \end{align}
for all $\eta$, $\xi\in\R^\dmn$ and almost every $x\in\R^\dmn\setminus\partial\Omega$.
\end{defn}

The main results (Theorems~\ref{thm:Dirichlet} and~\ref{thm:Dirichlet:Lp}) of the present paper hold in the full generality of real elliptic matrices, with no further restrictions on~$A$.

\begin{defn}[Weak solutions]\label{dfn:weak}
We will consider solutions to $-\Div A\nabla u=h$ in the case where $h$ is the divergence of some vector field~$\vec G$.

We will use the weak definition of $-\Div A\nabla u=-\Div \vec G$; that is, if $u:\Omega\to\R$ is a real-valued function, we say that $-\Div A\nabla u=-\Div \vec G$ in $\Omega$ in the weak sense if, for every $\varphi:\Omega\to \R$ that is smooth and supported in a compact subset of~$\Omega$, we have that
\begin{equation*}\int_\Omega\nabla\varphi\cdot A\nabla u=\int_\Omega\nabla\varphi\cdot \vec G.\end{equation*}
\end{defn}
Recall that $A$ is elliptic and therefore invertible.
For notational convenience we will usually take our data to be $\vec H=A^{-1}\vec G$; thus, we consider solutions to
$-\Div A\nabla u=-\Div (A\vec H)$, meaning functions $u$ that satisfy
\begin{equation}\label{eqn:weak}
\int_\Omega\nabla\varphi\cdot A\nabla u=\int_\Omega\nabla\varphi\cdot A\vec H
\quad\text{for all }\varphi\in C^\infty_0(\Omega).\end{equation}
This notational convention will considerably simplify the statement of Theorem~\ref{thm:Dirichlet} below.

We now state three known facts concerning solutions to $-\Div A\nabla u=0$ or $-\Div A\nabla u=-\Div(A\vec H)$ that we will reference in the statement of our main results.

First, for any domain $\Omega$ with $\bdmn$-Ahlfors regular boundary, $0<\bdmn\leq\dmnMinusOne$, and for any real elliptic operator $-\Div A\nabla $ as above, for every function $f$ defined on $\partial\Omega$ that is continuous and compactly supported, there is a unique solution $u$ to the Dirichlet problem
\begin{equation}\label{eqn:Dirichlet:cts}
\left\{\begin{gathered}\begin{aligned}
-\Div A\nabla u&=0&&\text{in }\Omega,\\
u&=f&&\text{on }\partial\Omega,\end{aligned}\\
u\text{ is continuous and bounded on $\overline\Omega$}
.\end{gathered}
\right.\end{equation}
That this holds for $\bdmn=\dmnMinusOne$ is the fact that Ahlfors regular boundaries satisfy the classic Wiener criterion \cite{Wie24}; that it holds for $\bdmn<\dmnMinusOne$ is \cite[Lemma~9.4]{DavFM21}.

The second result we wish to mention is Meyers's reverse H\"older estimate for gradients. This estimate states that, if $A$ satisfies the conditions~(\ref{eqn:elliptic:introduction}--\ref{eqn:elliptic:bounded:introduction}), then there is a $C$ and $\varepsilon>0$ depending only on the dimension $\dmn$ (and $\bdmn$) and on the ellipticity constants $\lambda$ and~$\Lambda$ such that, if $-\Div A\nabla u=-\Div (A\vec H)$ in $B(x,r)\subseteq \Omega$, and if $\vec H\in L^{2+\varepsilon}(B(x,r))$, then $\nabla u\in L^{2+\varepsilon}(B(x,r/2))$ and
\begin{equation}\label{eqn:Meyers}
\biggl(\fint_{B(x,r/2)} |\nabla u|^{2+\varepsilon}\biggr)^{1/(2+\varepsilon)}
\leq C\biggl(\fint_{B(x,r)} |\nabla u|^{2}\biggr)^{1/2}
+C\biggl(\fint_{B(x,r)} |\vec H|^{2+\varepsilon} \biggr)^{1/(2+\varepsilon)}
.\end{equation}
The same of course is true for solutions to $-\Div A^*\nabla u=-\Div (A^*\vec H)$, possibly with a different value of~$\varepsilon$.

Finally, we wish to state the boundary De Giorgi-Nash estimate. This estimate is well known in the case $\bdmn=\dmn-1$ and is still true if $0<\bdmn<\dmnMinusOne$; see \cite[Lemmas 8.13 and~8.16]{DavFM21}. It states that there exists an $\alpha>0$ and $C>0$ such that, if $-\Div A\nabla u=0$ in~$\Omega$, if $u=0$ on $\partial\Omega\cap B(\xi,r)$ for some $\xi\in\partial\Omega$ and some $r>0$, and if $x\in B(\xi,r/2)\cap\Omega$, then
\begin{equation}\label{eqn:bdry:DGN}
|u(x)|\leq \biggl(\frac{|x-\xi|}{r}\biggr)^\alpha
\frac{C}{r^{\bdmn+1}}
\int_{\Omega\cap B(\xi,r)} |u(y)|\, \dist(y,\partial\Omega)^{1+\bdmn-\pdmn}\,dy.\end{equation}

\section{A brief history of the Poisson-Dirichlet problem in spaces with fractional smoothness}\label{sec:history}

In this section we will briefly review the history of the problems which we consider in this note. That is, we will review the history of the Poisson-Dirichlet problem, with one degree of smoothness for the solutions, and with between zero and one degree of smoothness of the boundary values, for second order linear elliptic differential equations, with real coefficients and without lower order terms, in domains in Euclidean space.

\subsection{Higher codimension and energy solutions}\label{sec:energy}

Let $0<\bdmn<\dmn$, $\dmn\geq 2$, let $\Omega\subset\R^\dmn$ be such that $\partial\Omega$ is $\bdmn$-Ahlfors regular for some real number $\bdmn$ with $0<\bdmn<\dmn$, and let $A$ satisfy the ellipticity condition (\ref{eqn:elliptic:introduction}--\ref{eqn:elliptic:bounded:introduction}).
The Lax-Milgram lemma yields unique solvability of the Poisson problem
\begin{equation}
\label{eqn:Poisson}
-\Div A\nabla u_0=-\Div (A\vec H_0)\text{ in }\Omega,\quad u_0=0\text{ on }\partial\Omega,
\end{equation}
for all $\vec H_0\in L^2(\Omega,w)$, and with the estimate
\begin{equation}\label{eqn:Lax-Milgram:estimate}
\|\nabla u_0\|_{L^2(\Omega,w)}\leq C\|\vec H_0\|_{L^2(\Omega,w)}\end{equation}
where $u_0=0$ on $\partial\Omega$ in the sense that $u_0$ lies in $\dot W^{1,2}_0(\Omega;w)$, the closure of $C^\infty_0(\Omega)$ in ${\dot W^{1,2}(\Omega;w)}$ (the Sobolev space of functions with gradient in $L^{2}(\Omega;w)$), and where $w$ is the measure given by $dw=\dist(\,\cdot\,,\partial\Omega)^{\bdmn+1-\pdmn}\,dm$, $m$ denoting Lebesgue measure. This method of solution is extremely general, relying only on the fact that $\dot W^{1,2}_0(\Omega;w)$ is a Hilbert space and on the upper and lower bounds
\[\int_\Omega\nabla \varphi\cdot A\nabla \psi\leq \Lambda\|\varphi\|_{\dot W^{1,2}_0(\Omega;w)}\|\psi\|_{\dot W^{1,2}_0(\Omega;w)},
\qquad
\int_\Omega \nabla u\cdot A\nabla u \geq \lambda \|u\|_{\dot W^{1,2}_0(\Omega;w)}^2
\]
on functions $\varphi$, $\psi$, $u\in \dot W^{1,2}_0(\Omega;w)$. Indeed this construction is valid in many situations beyond that of this paper, such as systems and higher order equations.

The estimate
\begin{equation}\label{eqn:Lax-Milgram:Meyers:estimate}
\|\nabla u_0\|_{L^{p}(\Omega,w)}\leq C\|\vec H_0\|_{L^{p}(\Omega,w)}\end{equation}
is also true for all $p$ sufficiently close to~$2$. This can be derived by using the fact that $\dot W^{1,p}_0(\Omega;w)$, $1<p<\infty$, is an interpolation scale and using \v{S}ne\v{\i}berg's lemma.

The problem~\eqref{eqn:Dirichlet:Poisson:introduction}, with the notation of formula~\eqref{eqn:weak}, becomes
\begin{equation}\label{eqn:Dirichlet:Poisson}
-\Div (A\nabla u)=-\Div (A\vec H)\text{ in }\Omega,\quad
u=f\text{ on }\partial\Omega.\end{equation}
If $\partial\Omega$ is $\bdmn$-Ahlfors regular, $0<\bdmn<\dmn$, then we can extend an arbitrary $f\in \dot B^{p,p}_{1-1/p}(\partial\Omega)$ to a function $F$ in $\dot W^{1,p}(\Omega;w)$ of comparable norm. (The extension result is valid even if $\bdmn$ is not an integer.)  See \cite{JonW84} if $\bdmn=\dmnMinusOne$, \cite[Chapter~7]{DavFM21} and \cite[Chapter~8]{DavFM20p} if $p=2$, and \cite{BarMP25pB} in the general case.
Letting $F$ be this extension, $u_0$ the solution to the homogeneous Poisson problem~\eqref{eqn:Poisson} with data $\vec H_0=\vec H-\nabla F$, and $u=u_0+F$ yields a solution $u$ to the Poisson-Dirichlet problem~\eqref{eqn:Dirichlet:Poisson} with the estimate
\begin{equation}\label{eqn:energy:estimate}
\|\nabla u\|_{L^p(\Omega,w)}\leq C\|\vec H\|_{L^p(\Omega,w)}+C\|f\|_{\dot B^{p,p}_{1-1/p}(\partial\Omega)}.
\end{equation}
(Trace theorems often allow us to clarify the sense in which $u=f$ on~$\partial\Omega$; see the above references.)
This construction is very well known in the codimension 1 case ($\bdmn=\dmnMinusOne$) and is explained in detail in the general case in \cite[Lem\-ma~9.1]{DavFM21}.

We emphasize that up until the present work, in the higher codimension case, this construction of solutions to the Poisson-Dirichlet problem, restricted to boundary data in the particular Besov space $\dot B^{p,p}_{1-1/p}(\partial\Omega)$ for $p$ sufficiently close to~$2$, is the only known result for the Dirichlet-Poisson problem with boundary data of fractional smoothness. All of the following results in this section apply only to coefficients that satisfy the ellipticity conditions (\ref{eqn:elliptic:introduction}--\ref{eqn:elliptic:bounded:introduction}) in the particular case $\bdmn=\dmnMinusOne$.

\subsection{The Poisson problem with Sobolev estimates}\label{sec:Sobolev}

A natural line of study is to seek conditions on $\Omega$ and $A$ such that the estimate \eqref{eqn:Lax-Milgram:Meyers:estimate} holds for a wider range of~$p$.

The estimate $\|\nabla u_0\|_{L^p(\Omega)}\leq C\|\vec H_0\|_{L^p(\Omega)}$ for all $1<p<\infty$ has been established in the case where $A$ satisfies (\ref{eqn:elliptic:introduction}--\ref{eqn:elliptic:bounded:introduction}) (with $\bdmn=\dmnMinusOne$) and lies in the space $VMO$ of vanishing mean oscillation, and where $\partial\Omega$ is a bounded $C^1$ domain. See \cite{DiF96,KinZ01,AusQ02}. For any fixed $p\in (1,\infty)$, the result continues to hold if $A$ and $\Omega$ are allowed to deviate slightly from these assumptions; see \cite{ByuW04,Byu05} for the details.

\subsection{The Poisson-Dirichlet problem with weighted Sobolev estimates}\label{sec:MMS10}

Recall that in Section~\ref{sec:energy} we used the fact that a function  $f\in {\dot B^{p,p}_{1-1/p}(\partial\Omega)}$ may be extended to a function $F$ in the Sobolev space $\dot W^{1,p}(\R^\dmn)$.
It is also possible to extend functions in ${\dot B^{p,p}_{s}(\partial\Omega)}$ for a range of~$s$, and not merely the particular value $s=1-1/p$, if we take our extensions to lie in a \emph{weighted} Sobolev space with an appropriate weight. In particular, if $0<s<1$, $1\leq p\leq \infty$, and $\partial\Omega$ is $\pdmnMinusOne$-Ahlfors regular, then any $f\in \dot B^{p,p}_{s}(\partial\Omega)$ may be extended to a function $F$ with
\begin{equation*}\int_{\R^\dmn} |\nabla F(x)|^p\dist(x,\partial\Omega)^{p-1-ps}\,dx\approx \|f\|_{\dot B^{p,p}_{s}(\partial\Omega)}^p.\end{equation*}
See \cite{Usp61,Nik77B,Sha85,NikLM88,MitT06,Kim07,MazMS10,DavFM21,BarMP25pB}. The case of bounded Lipschitz domains is \cite[Proposition~4.1]{MitT06} and the $l=1$ case of \cite[Theorem~2.10]{Kim07}, but the result was known earlier in smoother domains.

This suggests consideration of the Dirichlet-Poisson problem~\eqref{eqn:Dirichlet:Poisson} with the estimate
\begin{multline}\label{eqn:MMS10:estimate}\int_\Omega |\nabla u(x)|^p\dist(x,\partial\Omega)^{p-1-ps}\,dx\\\leq C \int_\Omega |\vec H(x)|^p\dist(x,\partial\Omega)^{p-1-ps}\,dx+C\|f\|_{\dot B^{p,p}_{s}(\partial\Omega)}^p.\end{multline}
This problem was studied in \cite{MitT06} (which in fact considered the case where $L$ is the Hodge Laplacian in a Lipschitz domain in a compact Riemann manifold) and \cite{MazMS10} (which in fact studied the problem for elliptic systems of arbitrary even order, not only for the second-order equations that are the focus of the present note).

The paper \cite{MazMS10} established a result similar to those of the previous section. Specifically, if $A$ satisfies (\ref{eqn:elliptic:introduction}--\ref{eqn:elliptic:bounded:introduction}) (with $\bdmn=\dmnMinusOne$) and lies in $VMO$, and if $\Omega$ is a bounded Lipschitz domain whose unit outward normal also lies in~$VMO$ (this is a moderate generalization of the $C^1$ condition), then the estimate~\eqref{eqn:MMS10:estimate} holds for any $1<p<\infty$ and any $0<s<1$. Moreover, again as in the unweighted case of the previous section, for any fixed $p$ and~$s$ in this range, if $A$ is sufficiently close to~$VMO$ and $\Omega$ is a Lipschitz domain with sufficiently small constant, then the estimate~\eqref{eqn:MMS10:estimate} holds.

The $p=2$, $s=1/2$ case of these results is the energy solutions of Section~\ref{sec:energy}, while the $s=1-1/p$ case is the previous section.

The paper \cite{MitT06} established the estimate~\eqref{eqn:MMS10:estimate} in general bounded Lipschitz domains (without the $VMO$ or near-$VMO$ condition on the normals) but for a restricted range of $s$, $p$, and~$A$. Specifically, they required the coefficients $A$ to be symmetric and to satisfy the Dini-type condition
\begin{equation*}|A(x)-A(y)|\leq \omega(|x-y|), \qquad \int_0^1 \frac{\sqrt{\omega(t)}}{t}\,dt<\infty\end{equation*}
which is certainly stronger than lying in~$VMO$. They were then able to establish the estimate~\eqref{eqn:MMS10:estimate} for all $s$ and $p$ such that the point $(s,1/p)$ lies in the hexagonal region illustrated in Figure~\ref{fig:JK95} for some $\varepsilon$, $\varepsilon^*$, $\mathfrak{a}$, $\mathfrak{a}^*>0$.\footnote{The paper \cite{MitT06}, and its predecessors \cite{MitT00A,MitT05}, do not work with arbitrary second-order operators in Lipschitz domains in $\R^\dmn$, but with the Hodge Laplacian $\Delta_g$ defined on a compact Riemannian manifold with metric~$g$; more precisely, they work with the Schr\"odinger-type $L=-\Delta_g+V$ for a nonnegative potential~$V$.

We may return to the Euclidean situation of the present paper as follows. Given a bounded Lipschitz domain $\Omega$ in~$\R^\dmn$, one may easily construct a compact $\dmn$-dimensional manifold $M$ with $\Omega$ contained in one coordinate chart. In local coordinates the Hodge Laplacian is given by $\sqrt{|g|}\Delta_g=\Div(\sqrt{|g|} g^{-1}\nabla)$, and so we can attain any real symmetric coefficient matrix $A$ by choosing an appropriate metric on the chart containing~$\Omega$. The technical requirements on the potential $V$ in \cite{MitT00A,MitT05,MitT06} may be dealt with by choosing $V=0$ in a neighborhood of~$\Omega$ and $V>0$ in appropriate subsets of $M\setminus\Omega$.}

\def\figurealpha{0.4}
\def\figureeps{0.05}
\def\figuredimen{3}
\begin{figure}
\begin{tikzpicture}[scale=3]
\begin{scope}[shift={(0,0)}]
\draw [->] (-0.3,0)--(1.3,0) node [below ] {$\vphantom{1}s$};
\draw [->] (0,-0.3)--(0,1.3) node [ left] {$1/p$};
\fill  (1,1) -- (1,1/2-\figureeps) -- (\figurealpha,0) -- (0,0) -- (0,1/2+\figureeps) -- (1-\figurealpha,1) -- cycle;
\node [right] at (1,1/2-\figureeps) {$(1,\frac12-\varepsilon^*)$}; \node at (1,1/2-\figureeps) {$\circ$};
\node [right] at (1,1) {$(1,1)$}; \node at (1,1) {$\circ$};
\node [left] at (0,1/2+\figureeps) {$\frac12+\varepsilon$}; \node at (0,1/2+\figureeps) {$\circ$};
\node [below] at (\figurealpha,0) {$\phantom{1}\mathfrak{a}$}; \node at (\figurealpha,0) {$\circ$};
\node [above left] at (1-\figurealpha,1) {$(1-\mathfrak{a}^*,1)$}; \node at (1-\figurealpha,1) {$\circ$};
\end{scope}
\end{tikzpicture}
\caption{Values of the parameters $(s,1/p)$ such that the estimate~\eqref{eqn:MMS10:estimate} holds, for solutions to the problem~\eqref{eqn:Dirichlet:Poisson} in a Lipschitz domain, where $A$ satisfies a Dini-type continuity condition. This region includes the edges at $p=1$ and $p=\infty$ but does not include the edges at $s=0$ or $s=1$.}
\label{fig:JK95}
\end{figure}
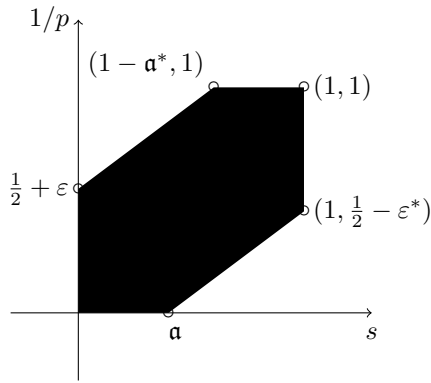

\subsection{General weights}

More general weights may be considered. By \cite{AdiMP21}, if $1<p<\infty$, if $\omega$ is a Muckenhoupt $A_p$ weight, if $A$ is symmetric and the $BMO$ norm of $A$ is small enough, and if $\Omega$ is a Lipschitz domain with sufficiently small Lipschitz constant, then solutions to the Poisson problem~\eqref{eqn:Poisson} satisfy the estimate
\begin{equation}\label{eqn:estimate:Lp:weight}
\|\nabla u_0\|_{L^p(\Omega;\omega)}\leq C
\|\vec H_0\|_{L^p(\Omega;\omega)}
.\end{equation}
The result was generalized to quasi-convex domains in \cite[Theorem~1.10]{YanYY22}. This builds on the earlier results \cite{Phu11,AdiP16}, among others, which proved the same result under stronger assumptions on~$\omega$.

As noted in \cite[formula~(2.5)]{BreM13} and the proof of~\cite[Lemma~2.3]{MitT06}, if $\Omega\subset\R^\dmn$ is a Lipschitz domain, then $\omega(x)=\dist(x,\partial\Omega)^{p-1-ps}$ is a Muckenhoupt $A_p$ weight for any $1<p<\infty$ and any $0<s<1$, and so the results of the previous section fall under this umbrella.

\subsection{Sharpness of norms}\label{sec:sharp}

The choice of norm on $f$ in Section~\ref{sec:MMS10} is sharp in the sense that, if the left-hand side is finite, then the boundary values of $u$ exist and lie in the Besov space $\dot B^{p,p}_{s}(\partial\Omega)$; indeed we have the reverse estimate
\begin{equation*}\|u\big\vert_{\partial\Omega}\|_{\dot B^{p,p}_{s}(\partial\Omega)}^p\leq C\int_\Omega |\nabla u(x)|^p\dist(x,\partial\Omega)^{p-1-ps}\,dx\end{equation*}
for \emph{any} function $u$ such that the right hand side is finite. This was proven in the generality of Lipschitz domains in \cite[Proposition~4.1]{MitT06} and as the $l=1$ case of \cite[Theorem~2.10]{Kim07}; similar results were already known in smoother domains (see \cite{Liz60,Usp61,Nik77B,Sha85,NikLM88}) or in the unweighted case $p-1-ps=0$ (see \cite{JonW84}).

The estimates on $\vec H$ (or $\vec H_0$) in all of the preceding sections are also sharp. Indeed, if $\Div(A\vec\Phi)=\Div(A\vec H_0)$, then the solution $u$ with data $\vec \Phi$ must equal the solution with data~$\vec H_0$, and so the estimate~\eqref{eqn:estimate:Lp:weight} may be strengthened to
\begin{equation*}\|\nabla u_0\|_{L^p(\Omega;\omega)}\leq C
\inf\bigl\{\|\vec\Phi\|_{L^p(\Omega;\omega)}:\Div (A\vec\Phi-A\vec H_0)=0\bigr\}.\end{equation*}
In particular, $\Div(A\vec H_0)=\Div (A\nabla\vec u)$, and so we may take $\vec\Phi=\nabla u_0$. Thus the reverse estimate
\begin{equation*}
\inf\bigl\{\|\vec\Phi\|_{L^p(\Omega;\omega)}:\Div (\vec\Phi-\vec H_0)=0\bigr\}
\leq C\|\nabla u_0\|_{L^p(\Omega;\omega)}
\end{equation*}
is trivially satisfied. The bound~\eqref{eqn:estimate:Lp:weight} thus implies an equivalence of norms
\begin{equation*}
\|\nabla u_0\|_{L^p(\Omega;\omega)}
\approx
\inf\bigl\{\|\vec\Phi\|_{L^p(\Omega;\omega)}:\Div (\vec\Phi-\vec H_0)=0\bigr\}.\end{equation*}

\subsection{Transversally constant coefficients and weighted averaged Sobolev estimates}\label{sec:t-independent}

Moving beyond Dini-type or VMO coefficients to another class of variable coefficient elliptic matrices, in \cite{BarM16A,AmeA18} the authors considered coefficients constant in the direction transverse to the boundary (but merely bounded measurable in all other directions). In \cite{BarM16A} (written by two of the authors of the present paper) we showed that, if $\Omega$ is the domain above a Lipschitz graph, and for such coefficients that in addition are real symmetric, we still have well posedness of the Poisson-Dirichlet problem~\eqref{eqn:Dirichlet:Poisson}, with the estimate \eqref{eqn:MMS10:estimate} replaced by the estimate
\begin{multline}\label{eqn:BarM16A:estimate}\int_\Omega
\biggl(\fint_{B(x,\dist(x,\partial\Omega)/2)}|\nabla u|^\beta\biggr)^{p/\beta} \dist(x,\partial\Omega)^{p-1-ps}\,dx\\\leq C \int_\Omega \biggl(\fint_{B(x,\dist(x,\partial\Omega)/2)}|\vec H|^\beta\biggr)^{p/\beta} \dist(x,\partial\Omega)^{p-1-ps}\,dx+C\|f\|_{\dot B^{p,p}_{s}(\partial\Omega)}^p,\end{multline}
for $\beta$ sufficiently close to~$2$,
and for $(s,1/p)$ in the hexagon in Figure~\ref{fig:MM04}. (In the case where $A$ is real, constant in the transverse direction, but not necessarily symmetric, we have well posedness for a smaller range; namely, the range of Figure~\ref{fig:MM04} for some $\varepsilon$, $\varepsilon^*$ greater than~$-1/2$. See the region at the bottom of Figure~\ref{fig:extrapolation} below.)

\def\figurealpha{0.4}
\def\figureeps{0.05}
\def\figuredimen{3}
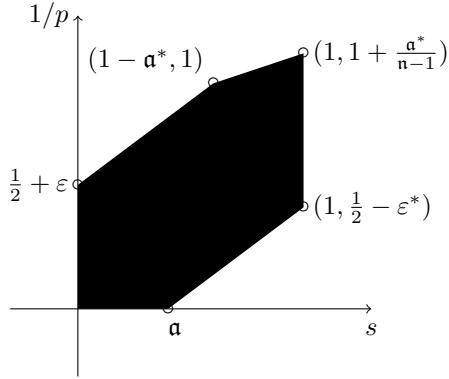
\begin{figure}
\begin{tikzpicture}[scale=3]
\begin{scope}[shift={(0,0)}]
\draw [->] (-0.3,0)--(1.3,0) node [below ] {$\vphantom{1}s$};
\draw [->] (0,-0.3)--(0,1.3) node [ left] {$1/p$};
\fill  (1,1) -- (1,1/2-\figureeps) -- (\figurealpha,0) -- (0,0) -- (0,1/2+\figureeps) -- (1-\figurealpha,1) --(1,1+\figurealpha/\figuredimen) -- cycle;
\node [right] at (1,1/2-\figureeps) {$(1,\frac12-\varepsilon^*)$}; \node at (1,1/2-\figureeps) {$\circ$};
\node [right] at (1,1+\figurealpha/\figuredimen) {$(1,1+\frac{\mathfrak{a}^*}{\dmnMinusOne})$}; \node at (1,1+\figurealpha/\figuredimen) {$\circ$};
\node [left] at (0,1/2+\figureeps) {$\frac12+\varepsilon$}; \node at (0,1/2+\figureeps) {$\circ$};
\node [below] at (\figurealpha,0) {$\phantom{1}\mathfrak{a}$}; \node at (\figurealpha,0) {$\circ$};
\node [above left] at (1-\figurealpha,1) {$(1-\mathfrak{a}^*,1)$}; \node at (1-\figurealpha,1) {$\circ$};
\end{scope}
\end{tikzpicture}
\caption{Values of the parameters $(s,1/p)$ such that solutions to Dirichlet-Poisson problem~\eqref{eqn:Dirichlet:Poisson} in the domain above a Lipschitz graph, for real symmetric coefficients constant in the vertical direction, satisfy the estimate~\eqref{eqn:BarM16A:estimate}.}
\label{fig:MM04}
\end{figure}

The $\beta$-power averages were inspired by the modified nontangential maximal function of \cite{KenP93} and were introduced for much the same reason: if $A$ lacks smoothness properties, then solutions $u$ to even the homogeneous problem $\Div A\nabla u=0$ no longer have gradients with high degrees of integrability. The gradient $\nabla u$ may be locally in $L^\beta$ only for $\beta$ close to~$2$, and so these weighted averaged spaces are more appropriate to such differential equations. Even for $p$ much larger than~$\beta$, these spaces preserve the trace result
\begin{equation*}\|f\|_{\dot B^{p,p}_{s}(\partial\Omega)}^p\leq C\int_\Omega
\biggl(\fint_{B(x,\dist(x,\partial\Omega)/2)}|\nabla u|^\beta\biggr)^{p/\beta} \dist(x,\partial\Omega)^{p-1-ps}\,dx\end{equation*}
and so the estimate~\eqref{eqn:BarM16A:estimate} is sharp.
See \cite[Theorem~6.3]{BarM16A}.

As mentioned above, we used layer potential techniques, which introduce harsh restrictions both on the coefficients and on the domains. Our results also used crucially the endpoint results (which may be viewed as the $s=0$ and $s=1$ boundaries of the hexagon in Figure~\ref{fig:MM04}); that is, with the techniques of \cite{BarM16A}, the results for Besov boundary values do not stand alone, but require the results for boundary data in $L^p(\partial\Omega)$ (or in the Sobolev space $\dot W^{1,p}(\partial\Omega)$, that is, at the smoothness~$1$ endpoint) and require the layer potential bounds.

The paper \cite{AmeA18} applies the so-called `first-order approach' (developed in the end point theory of $L^p(\partial\Omega)$ boundary data) to the Dirichlet problem
\begin{equation}\label{eqn:Dirichlet}
-\Div A\nabla u=0\text{ in }\Omega,\quad
u=f\text{ on }\partial\Omega\end{equation}
with $f$ in a fractional smoothness space, either the Besov spaces considered in this note or the related Hardy-Sobolev spaces.

Like \cite{BarM16A}, they considered coefficients constant in the direction transverse to the boundary. Again like \cite{BarM16A}, their direct well posedness results for boundary data of fractional smoothness require end point results, that is, results for boundary data in $L^p(\partial\Omega)$ or $\dot W^{1,p}(\partial\Omega)$.

However, unlike \cite{BarM16A}, the methods of \cite{AmeA18} encompass general complex coefficients, and even elliptic systems rather than equations. \cite{BarM16A} is confined to equations, and also to a restricted class of coefficients (specifically, those that satisfy the De Giorgi-Nash condition; this includes all real coefficients but misses many complex coefficients). \cite{AmeA18} thus can build upon endpoint results for a much broader class of coefficients, including complex self-adjoint coefficients, and constant or triangular coefficients in the half-space $\R^\dmn_+$, to produce results for Besov data.

\subsection{Non-sharp estimates}\label{sec:nonsharp}

We mention one further class of results concerning the Poisson problem with rough variable coefficients. The Poisson problem \eqref{eqn:Poisson} with the estimate
\begin{equation*}\int_\Omega |\nabla u_0(x)|^p\dist(x,\partial\Omega)^{p-1-ps}\,dx\leq C \int_\Omega |\vec H_0(x)|^p\dist(x,\partial\Omega)^{p-1-ps}\,dx\end{equation*}
has sharp estimates, is well posed in many circumstances, and is well-adapted to be of use in analyzing the Dirichlet problem with boundary data in the Besov space~$\dot B^{p,p}_{s}(\partial\Omega)$.

There is no corresponding estimate for the Poisson problem that is well adapted to the $L^p(\partial\Omega)$-Dirichlet problem, is valid in many (or even any known) circumstances, and is sharp in the sense of Section~\ref{sec:sharp}.

However, it is possible to find estimates for the Poisson problem, adapted to the $L^p(\partial\Omega)$-Dirichlet problem, that do hold in many circumstances; these estimates are simply not sharp. Such estimates were found in the works
\cite{HofMayMou15,Bar20p,MouPT22p}, either implicitly as part of the method of attack on the $L^p$-Dirichlet problem \cite{HofMayMou15,Bar20p} or explicitly as a main result \cite{MouPT22p}.

Specifically, these papers established the estimates
\begin{equation*}\|\widetilde Nu\|_{L^p(\partial\Omega)}\leq C\|\mathcal{C}_2(\vec H)\|_{L^p(\partial\Omega)}\end{equation*}
and
\begin{equation*}\|\widetilde N(\nabla u)\|_{L^p(\partial\Omega)}\leq C\|\mathcal{C}_2(\vec H/\dist(\,\cdot\,,\partial\Omega))\|_{L^p(\partial\Omega)}\end{equation*}
where $\widetilde N$ and $\mathcal{C}_2$ are the modified nontangential maximal function and tent space operator common in the $L^p$ theory. (We refer the reader to \cite{HofMayMou15,Bar20p,MouPT22p} for precise definitions.)

These estimates do establish well posedness of a form of the Poisson problem, but they are not sharp in the sense that there always exist~$u$ such that\footnote{For example, in the case~\eqref{eqn:not:Dirichlet}, we may take $u_k(x)=\eta(x)\sin(kx_1)$, where $\eta$ is a nonnegative smooth cutoff function. Then $\widetilde Nu_k\leq \widetilde N\eta$ and so is uniformly bounded in~$L^p$, but $\nabla u_k$ grows without bound as $k\to\infty$ and so the left hand term $\mathcal{C}_2(\nabla u_k)$ cannot be bounded uniformly in~$k$.

In the case~\eqref{eqn:not:regularity}, we may choose $u$ such that $\nabla u$ is a nonzero constant in a ball centered at the boundary. Then $\mathcal{C}_2(\nabla u/\dist(\,\cdot\,,\partial\Omega))$ involves an integral of $1/\dist(\,\cdot\,,\partial\Omega)$ and so will diverge in a nontrivial set, while $u$ may be cut off outside that ball in such a way that $\widetilde N(\nabla u)$ is bounded and compactly supported.}
\begin{equation}\label{eqn:not:Dirichlet}\|\mathcal{C}_2(\nabla u)\|_{L^p(\partial\Omega)} \not\leq C\|\widetilde Nu\|_{L^p(\partial\Omega)}\end{equation}
or
\begin{equation}\label{eqn:not:regularity}\|\mathcal{C}_2(\nabla u/\dist(\,\cdot\,,\partial\Omega))\|_{L^p(\partial\Omega)} \not\leq C\|\widetilde N(\nabla u)\|_{L^p(\partial\Omega)}.\end{equation}

\subsection{The Laplace operator and other estimates}\label{sec:Laplace}

Thus far we have focused on estimates of solutions $u$ to the Poisson problem that satisfy local Sobolev estimates, that is, solutions where the gradient $\nabla u$ satisfies an estimate that locally reduces to membership in $L^2$ or~$L^p$. We have done this because the main new results of this note (to be described in Sections \ref{sec:Poisson} and~\ref{sec:extrapolation}) are phrased in these terms.

However, many important results, particularly the earliest results, instead established estimates on solutions to the Poisson problem~\eqref{eqn:Dirichlet:Poisson:introduction} of the form
\begin{align}
\label{eqn:MM04:estimate}
\|u\|_{\dot B^{p,q}_{s+1/p}(\Omega)}&\leq C \|h\|_{\dot B^{p,q}_{s+1/p-2}(\Omega)}
+ C\|f\|_{\dot B^{p,q}_s(\partial\Omega)} \qquad\text{ or}
\\
\label{eqn:JK95:estimate}
\|u\|_{\dot F^{p,q}_{s+1/p}(\Omega)}&\leq C \|h\|_{\dot F^{p,q}_{s+1/p-2}(\Omega)}
+ C\|f\|_{\dot B^{p,p}_s(\partial\Omega)}
\end{align}
where $\dot B^{p,q}_\theta$ and $\dot F^{p,q}_\theta$ denote the Besov and Triebel-Lizorkin spaces.
(The precise statement of the estimates varies from paper to paper.)

It is often possible to pass from Besov or Triebel-Lizorkin estimates on solutions to weighted Sobolev estimates such as~\eqref{eqn:MMS10:estimate}; see, for example, \cite[Theorems 4.1 and~4.2]{JerK95}.

Estimates of the form \eqref{eqn:MM04:estimate} or~\eqref{eqn:JK95:estimate} are unsuited for solutions to $-\Div A\nabla u=-\Div (A\vec H)$, or even to $-\Div A\nabla u=0$, for rough coefficients~$A$. Such solutions $u$ do not generally display any higher regularity, and so cannot be expected to lie in $\dot B^{p,q}_{s+1/p}(\Omega)$ for $s+1/p>1$. See the discussion in \cite[Chapter~10]{BarM16A}, or compare the two papers \cite{MitT05} and \cite{MitT06}; both consider the Dirichlet-Poisson problem in Lipschitz domains for coefficients $A$ satisfying the same Dini-type condition, but \cite{MitT05} used the traditional estimate~\eqref{eqn:MM04:estimate} while \cite{MitT06} used the estimate~\eqref{eqn:MMS10:estimate}, and so \cite{MitT05} has a somewhat artifical upper bound on $s+1/p$ which \cite{MitT06} does not need to impose.

However, many important results are known in the context of the estimate \eqref{eqn:MM04:estimate} or~\eqref{eqn:JK95:estimate}.

The paper \cite{JerK95}, concerning only the Laplacian, is one of the foundational papers in the theory of the Poisson problem and of the Dirichlet problem with boundary data of fractional smoothness; indeed it preceded all of the results described in  Sections \ref{sec:Sobolev}--\nobreak\ref{sec:nonsharp}. That paper established the estimate~\eqref{eqn:JK95:estimate} in Lipschitz domains for $q=2$ and for $s$ and $p$ as in Figure~\ref{fig:JK95}. They provided examples establishing sharpness of these results (that is, of Lipschitz domains such that the estimate~\eqref{eqn:JK95:estimate} fails for any $(s,1/p)$ outside of the given region).
In the case of $C^1$ domains, they established this estimate for all $s\in (0,1)$ and all $p\in (1,\infty)$. (The $VMO$ condition of \cite{MazMS10} described in Section~\ref{sec:MMS10} is a mild generalization of the $C^1$ condition.)

The paper \cite{MayMit04A} removes the restriction on $q$ in the estimate \eqref{eqn:JK95:estimate} (and \eqref{eqn:MM04:estimate}) and expands the available values of $p$ and $s$ into the range $p<1$, specifically to the range $1/p<1+\frac{s+\mathfrak{a}^*-1}{\dmnMinusOne}$ (in Lipschitz domains) or $0\leq 1/p<1+\frac{s}{\dmnMinusOne}$ (in $C^1$ domains). See Figure~\ref{fig:MM04}. We observe that \cite{MayMit04A} precedes the papers \cite{BarM16A,AmeA18} described in Section~\ref{sec:t-independent} and that the paper \cite{BarM16A} explicitly works with the same ranges of $s$ and $p$ as the Lipschitz case of \cite{MayMit04A}.

The results of \cite{JerK95,MayMit04A} for $C^1$ domains with the estimate \eqref{eqn:MM04:estimate} or~\eqref{eqn:JK95:estimate} also hold in bounded convex domains; see \cite{MitMY10}.

\subsection{Hölder continuous boundary data}\label{sec:Holder}

We conclude this section with a discussion of the recent result \cite{CaoHMPZ23p}. This result considers the Dirichlet problem~\eqref{eqn:Dirichlet:cts} with H\"older continuous boundary data; this is the $p=\infty$ (bottom) edge of the region in Figure~\ref{fig:MM04}, and thus falls within the ambit of boundary value problems with data of fractional smoothness.

We remark that \cite{CaoHMPZ23p} is not the only paper to treat H\"older continuous boundary values, but in the setting of real elliptic coefficients, it contains (provably) the most general possible result. In particular, they place no assumptions on the coefficients beyond the ellipticity conditions (\ref{eqn:elliptic:introduction}--\ref{eqn:elliptic:bounded:introduction}) (with $\bdmn=\dmnMinusOne$). They showed that, if $\partial\Omega$ satisfies a criterion of Wiener type, then there is a $\varepsilon>0$ such that the solution $u$ to the continuous Dirichlet problem~\eqref{eqn:Dirichlet:cts} satisfies the estimate \begin{equation}\label{eqn:Holder}
\|u\|_{\dot C^s(\Omega)}\leq C\|f\|_{\dot C^s(\partial\Omega)}\end{equation}
for all $s>0$ sufficiently small, where $\dot C^s$ denotes spaces of H\"older continuous functions. Their result is sharp; they show that their Wiener-type criterion is a necessary, as well as a sufficient, condition for the estimate~\eqref{eqn:Holder} for arbitrary $f\in \dot C^s(\partial\Omega)$.

We remark that Ahlfors regularity of $\partial\Omega$ for $\bdmn=\dmnMinusOne$ (or indeed for any $\bdmn\in (\dmn-2,\dmn)$) implies their Wiener type criterion, and thus is a sufficient condition for H\"older well posedness.

The $p=\infty$, $\vec H=0$ analogue of the estimate~\eqref{eqn:BarM16A:estimate} is
\begin{equation*}\sup_{x\in\Omega } \biggl(\fint_{B(x,\dist(x,\partial\Omega)/2)}|\nabla u|^\beta\biggr)^{1/\beta} \dist(x,\partial\Omega)^{1-s}
\leq C\|f\|_{\dot C^s(\partial\Omega)}\end{equation*}
and it follows immediately (for $\beta$ sufficiently close to~$2$) from the estimate~\eqref{eqn:Holder} and the Caccioppoli and Meyers estimates for solutions to $-\Div A\nabla u=0$.

\section{The Poisson problem for very general domains and coefficients}\label{sec:Poisson}

Our first main theorem establishes well posedness of the Poisson-Dirichlet problem, for \emph{any} second order elliptic differential equation with real coefficients satisfying suitable upper and lower bounds, and with very mild conditions on the underlying domains.

\begin{thm}\label{thm:Dirichlet}
Let $\Omega\subseteq\R^\dmn$, $\dmn\geq 3$, be a connected open set whose boundary $\partial\Omega$ is $\pbdmn$-Ahlfors regular in the sense of Definition~\ref{dfn:adr} for some real number $\bdmn$ with $0<\bdmn\leq \dmn-1$. If $\bdmn=\dmnMinusOne$ then we additionally require that $\Omega$ satisfies the weak local John condition and the interior corkscrew condition.\footnote{If $\bdmn<\dmnMinusOne$ then, as noted in Remark~\ref{rmk:higher:John:corkscrew}, the weak local John and interior corkscrew conditions are automatically satisfied.}

We require that $\R^\dmn\setminus\Omega$ be unbounded, that is, either $\Omega$ is bounded\footnote{It is straightforward to establish that if $\Omega\subset\R^\dmn$ is bounded, then $\partial\Omega$ has dimension at least $\dmnMinusOne$; in our circumstances this means that $\bdmn=\dmnMinusOne$.} or both $\Omega$ and $\partial\Omega$ are unbounded.

Let $L=-\Div A\nabla$ be a linear differential operator associated to real, not necessarily symmetric coefficients $A$ that satisfy the positive definiteness and boundedness conditions (\ref{eqn:elliptic:introduction}--\ref{eqn:elliptic:bounded:introduction}).

There are then numbers $\mathfrak{a}$, $\mathfrak{a}^*\in (0,1]$ and $\varepsilon$, $\delta>0$ such that, if $(s,1/p)$ lies in the open pentagon indicated in Figure~\ref{fig:arbitrary}, and if $2-\delta<\beta<2+\varepsilon$, then
for all $f$ in the Besov space $\dot B^{p,p}_{s}(\partial\Omega)$ defined in Definition~\ref{dfn:Besov} above, and all
$\vec H\in L^\beta_{loc}(\Omega)$ such that the right hand side of the bound~\eqref{eqn:Besov:estimate} is finite, there is a unique solution $u$ to the Poisson-Dirichlet problem
\begin{equation}\label{eqn:Dirichlet:besov}
\left\{\begin{aligned}
-\Div A\nabla u&=-\Div(A\vec H)&&\text{in }\Omega\\
u&=f&&\text{on }\partial\Omega
\end{aligned}\right.\end{equation}
that satisfies
\begin{multline}\label{eqn:Besov:estimate}
\int_\Omega
\biggl(\fint_{B(x,\dist(x,\partial\Omega)/2)} |\nabla u|^\beta \biggr)^{p/\beta}
\dist(x,\partial\Omega)^ {\bdmn-\pdmn+p-ps}
\,dx
\\\leq
C\|f\|_{\dot B^{p,p}_{s}(\partial\Omega)}^p+
C\int_\Omega
\biggl(\fint_{B(x,\dist(x,\partial\Omega)/2)}|\vec H|^\beta \biggr)^{p/\beta}
\dist(x,\partial\Omega)^ {\bdmn-\pdmn+p-ps}
\,dx
\end{multline}
for some constant $C$ depending on $p$, $s$, $A$, and~$\Omega$, but not on $f$ or~$\vec H$. If $p=\infty$ and $0<s<\mathfrak{a}$, then there is a unique solution to the problem~\eqref{eqn:Dirichlet:besov} that satisfies
\begin{multline}\label{eqn:Besov:estimate:infinity}
\sup_{x\in\Omega}
\biggl(\fint_{B(x,\dist(x,\partial\Omega)/2)}|\nabla u|^\beta\biggr)^{1/\beta}
\dist(x,\partial\Omega)^ {1-s}
\\\leq
C\|f\|_{\dot B^{\infty,\infty}_{s}(\partial\Omega)}+
C\sup_{x\in\Omega}
\biggl(\fint_{B(x,\dist(x,\partial\Omega)/2)}|\vec H|^\beta\biggr)^{1/\beta}
\dist(x,\partial\Omega)^ {1-s}
.\end{multline}

Furthermore, we have that $\mathfrak{a}\geq \max(\alpha,1-\bdmn)$ and $\mathfrak{a}^*\geq \max(\alpha^*,1-\bdmn)$, where $\alpha$ and $\alpha^*$ are the numbers in the boundary De Giorgi-Nash estimate~\eqref{eqn:bdry:DGN}.

Finally, $2+\varepsilon$ may be taken to be the number in Meyers's reverse H\"older estimate for~$L$, and $2-\delta$ may be taken to be the H\"older conjugate of the number in Meyers's reverse H\"older estimate for~$L^*$.
\end{thm}

\def\figurealpha{0.3}
\def\figuredimen{2}
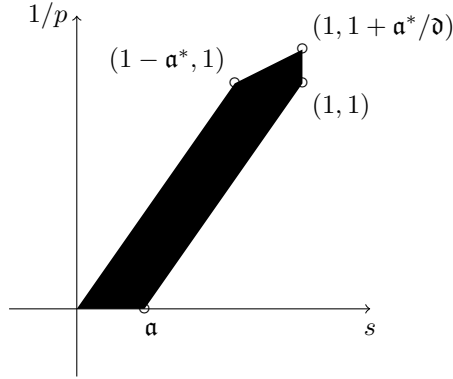
\begin{figure}
\begin{tikzpicture}[scale=3]
\begin{scope}[shift={(0,0)}]
\draw [->] (-0.3,0)--(1.3,0) node [below ] {$\vphantom{1}s$};
\draw [->] (0,-0.3)--(0,1.3) node [ left] {$1/p$};
\fill  (1,1) -- (\figurealpha,0) -- (0,0) -- (1-\figurealpha,1) --(1,1+\figurealpha/\figuredimen) -- cycle;
\node at (1,1) {$\circ$};
\node [below right] at (1,1) {$(1,1)$};
\node at (1,1+\figurealpha/\figuredimen) {$\circ$};
\node [above right] at (1,1+\figurealpha/\figuredimen) {$(1,1+\mathfrak{a}^*/\bdmn)$};
\node at (1-\figurealpha,1) {$\circ$};
\node [above left] at (1-\figurealpha,1) {$(1-\mathfrak{a}^*,1)$};
\node [below] at (\figurealpha,0) {$\phantom{1}\mathfrak{a}$}; \node at (\figurealpha,0) {$\circ$};
\end{scope}
\end{tikzpicture}
\caption{Values of the parameters $(s,1/p)$ such that the $\dot B^{p,p}_{s}(\partial\Omega)$-Dirichlet problem and associated Poisson problem for $L$ is solvable. By Remark~\ref{rmk:convex}, the illustrated region is convex.}
\label{fig:arbitrary}
\end{figure}

Expressed algebraically, the pair $(s,1/p)$ satisfies the conditions of Theorem~\ref{thm:Dirichlet} if either
\begin{align}
\label{eqn:thm:Dirichlet:quasi}
\frac{\bdmn}{\bdmn+\mathfrak{a}^*}&<p\leq 1&\text{and}&&\frac{\bdmn}{p}+1-\bdmn-\mathfrak{a}^*&<s<1  \\
\noalign{\noindent\text{or}}
\label{eqn:thm:Dirichlet:Banach}
1&\leq p\leq \infty &\text{and}&&\frac{1-\mathfrak{a}^*}{p}&<s<\mathfrak{a}+\frac{1-\mathfrak{a}}{p}.
\end{align}

We observe that, with the exception of the Sobolev estimate~\eqref{eqn:energy:estimate} for $p$ sufficiently close to~$2$ in Section~\ref{sec:energy}, and the Dirichlet problem with H\"older continuous boundary data in Section~\ref{sec:Holder}, all of the results for the Poisson and Dirichlet problems described in Section~\ref{sec:history} were shown exclusively in Lipschitz domains with boundary of codimension~$1$, or in domains with a higher degree of smoothness, and for coefficients satisfying some regularity assumption. Our result is the first result to produce well posedness for a broad range of $(s,1/p)$ for arbitrary real coefficients and for domains in the above generality.

Because fundamental solutions and Green's functions display different behavior in two dimensions and in higher dimensions, certain estimates in our proof of Theorem~\ref{thm:Dirichlet} work only if the ambient dimension $\dmn$ is at least three. The case $\dmn=2$ is open.

\begin{rmk}\label{rmk:convex}
Interpolation methods can generally be used to show that, if the $\dot B^{p,p}_{s}$-Dirichlet problem is well posed for $(s,1/p)=(s_0,1/p_0)$ and $(s,1/p)=(s_1,1/p_1)$, then it is well posed whenever $(s,1/p)$ lies on the line segment connecting the two points $(s_0,1/p_0)$ and $(s_1,1/p_1)$. Indeed such techniques were used in \cite{BarMP25pA} to prove Theorem~\ref{thm:Dirichlet}. Thus, regions such as those in Figure~\ref{fig:arbitrary} (or Figures~\ref{fig:JK95}--\ref{fig:MM04}) are generally expected to be convex.

The region of Figure~\ref{fig:arbitrary} is indeed convex. To see this, observe that the given region is a polygon, and thus we need only establish convexity in a neighborhood of each vertex.
The vertices at the origin, $(\mathfrak{a},0)$, $(1,1)$, and $(1,1+\mathfrak{a}^*/\bdmn)$ clearly have convex neighborhoods. We are left with the vertex at the point $(1-\mathfrak{a}^*,1)$.

The uppermost edge (connecting the points $(1-\mathfrak{a}^*,1)$ and $(1,1+\mathfrak{a}^*/\bdmn)$) has slope $1/\bdmn$. The adjacent edge connecting the point $(1-\mathfrak{a}^*,1)$ to the origin has slope $1/(1-\mathfrak{a}^*)$. But because $\mathfrak{a}^*\geq 1-\bdmn$, we see that $1/(1-\mathfrak{a}^*)\geq 1/\bdmn$ and so there must be a convex neighborhood of the vertex $(1-\mathfrak{a}^*,1)$ as well.
\end{rmk}

\section{An extrapolation result}\label{sec:extrapolation}

The result of the previous chapter is valid in the greatest known generality of coefficients and domains. It is essentially a consequence of the De Giorgi-Nash-Moser and Meyers inequalities which are a part of fundamental elliptic theory. As a consequence, we observe that the pentagon in Figure~\ref{fig:arbitrary} is necessarily smaller than the hexagon in Figure~\ref{fig:MM04}.

The goal of our project, however, was to exploit new well-posedness results for the Dirichlet problem in $L^p$ to obtain an enhanced range of Besov spaces. We now return to this task, and the theorem below will generalize and recover most of the results in Section~\ref{sec:history}, in particular, those in
Figure~\ref{fig:MM04}, and will go far beyond that.

To this end, recall that the Dirichlet problem is well-posed in $L^p$ for some $1<p<\infty$,  if for any  $f\in L^p(\partial\Omega)$ the solution $u$ to the continuous Dirichlet problem~\eqref{eqn:Dirichlet:cts} satisfies the nontangential estimate
\begin{equation}
\label{eqn:N:bound}
\|Nu\|_{L^q(\partial\Omega)}\leq C_q\|f\|_{L^q(\partial\Omega)}, \end{equation}
and if such a solution is unique. Here, $N$ is the nontangential maximal function given by
\begin{equation}\label{dfn:N}NF(x)=\sup \{|F(y)|:y\in\Omega,\>|x-y|<2\dist(y,\partial\Omega)\}
.\end{equation}

Well-posedness of the Dirichlet boundary value problem in $L^p$ is a powerful feature. It has long been known that it is equivalent to (quantitative) absolute continuity of harmonic measure with respect to the Hausdorff measure on the boundary; see \cite[Section~4]{HofL18} or \cite{Hof19} for a detailed discussion.
Under suitable connectivity assumptions and assumptions on the coefficients, well posedness of the $L^p$ Dirichlet problem for some $p$ is equivalent to uniform rectifiability for $\bdmn=n-1$; see \cite{AzzHMMT20,HofMMTZ21}.
We will highlight some recent work in this context in Section~\ref{sec:Lp}. For now, let us state our second main result.

\def\figureq{0.2}
\begin{figure}
\begin{tikzpicture}[scale=3]
\begin{scope}[shift={(0,0)}]
\draw [->] (-0.3,0)--(1.3,0);
\draw [->] (0,-0.3)--(0,1.3);
\node at (1,1) {$\circ$};
\node [below right] at (1,1) {$(1,1)$};
\node at (1,1+\figurealpha/\figuredimen) {$\circ$};
\node [above right] at (1,1+\figurealpha/\figuredimen) {$(1,1+\mathfrak{b}^*/\bdmn)$};
\node at (1-\figurealpha,1) {$\circ$};
\node [above left] at (1-\figurealpha,1) {$(1-\mathfrak{b}^*,1)$};
\node [left] at (0,\figureq) {$1/q$}; \node at (0,\figureq) {$\circ$};
\node [below] at (\figurealpha,0) {$\phantom{1}\mathfrak{a}$}; \node at (\figurealpha,0) {$\circ$};
\fill
(0,0) --(0,\figureq) --(1-\figurealpha,1) --(1,1+\figurealpha/\figuredimen) -- (1,1)
--(\figurealpha,0) --cycle;
\end{scope}
\begin{scope}[shift={(2,0)}]
\draw [->] (-0.3,0)--(1.3,0);
\draw [->] (0,-0.3)--(0,1.3);
\node at (1,1-\figureq) {$\circ$};
\node [right] at (1,1-\figureq) {$(1,1-1/q^*)$};
\node at (1,1+\figurealpha/\figuredimen) {$\circ$};
\node [above right] at (1,1+\figurealpha/\figuredimen) {$(1,1+\mathfrak{a}^*/\bdmn)$};
\node at (1-\figurealpha,1) {$\circ$};
\node [above left] at (1-\figurealpha,1) {$(1-\mathfrak{a}^*,1)$};
\node [below] at (\figurealpha,0) {$\phantom{1}\mathfrak{b}$}; \node at (\figurealpha,0) {$\circ$};
\fill
(1,1) --(1,1-\figureq) --(\figurealpha,0)--(0,0)
--(1-\figurealpha,1) --(1,1+\figurealpha/\figuredimen) --cycle;
\end{scope}
\end{tikzpicture}\\
\begin{tikzpicture}[scale=3]
\draw [->] (-0.3,0)--(1.3,0);
\draw [->] (0,-0.3)--(0,1.3);
\node at (1,1-\figureq) {$\circ$};
\node [right] at (1,1-\figureq) {$(1,1-1/q^*)$};
\node at (1,1+\figurealpha/\figuredimen) {$\circ$};
\node [above right] at (1,1+\figurealpha/\figuredimen) {$(1,1+\mathfrak{b}^*/\bdmn)$};
\node at (1-\figurealpha,1) {$\circ$};
\node [above left] at (1-\figurealpha,1) {$(1-\mathfrak{b}^*,1)$};
\node [left] at (0,\figureq) {$1/q$}; \node at (0,\figureq) {$\circ$};
\node [below] at (\figurealpha,0) {$\phantom{1}\mathfrak{b}$}; \node at (\figurealpha,0) {$\circ$};
\fill
(1,1) --(1,1-\figureq) --(\figurealpha,0) --(0,0)
--(0,\figureq) --(1-\figurealpha,1) --(1,1+\figurealpha/\figuredimen) --cycle;
\end{tikzpicture}
\caption{Values of the parameters $(s,1/p)$ such that the $\dot B^{p,p}_{s}(\partial\Omega)$-Dirichlet problem for $L$ is solvable, given solvability of the $L^q$-Dirichlet problem for $-\Div A\nabla$ (on the left) or the $L^{q^*}$-Dirichlet problem for $-\Div A^*\nabla$ (on the right) with nontangential estimates.  Because $\mathfrak{b}^*\geq 1-\bdmn+\bdmn/q$, by an argument similar to that in Remark~\ref{rmk:convex}, the illustrated region is convex.
In the common case that $q$ and $q^*$ both exist, the $\dot B^{p,p}_{s}(\partial\Omega)$-Dirichlet problem for $-\Div A\nabla$ is solvable for values of $(s,1/p)$ as in the bottom picture.
}
\label{fig:extrapolation}
\end{figure}
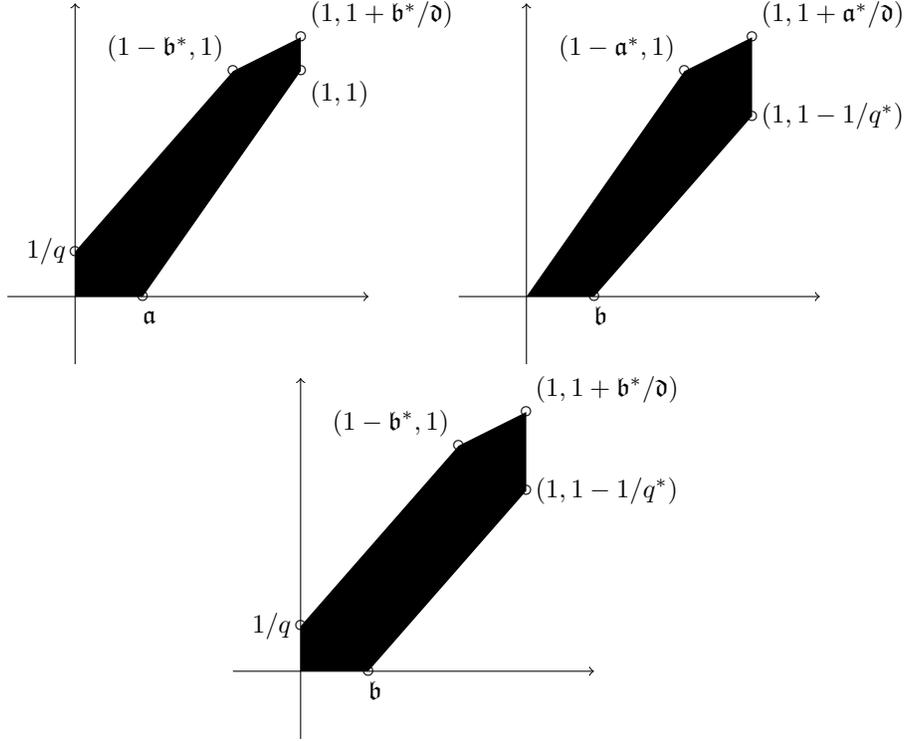

\begin{thm}\label{thm:Dirichlet:Lp}
Suppose that the conditions of Theorem~\ref{thm:Dirichlet} are satisfied.

Suppose furthermore that there is a $q\in (1,\infty)$ such that the $L^q$-Dirichlet problem is well posed. That is, suppose that there is a $q\in (1,\infty)$ and a $C_q<\infty$ such that if $f$ is continuous and compactly supported, then the solution $u$ to the problem~\eqref{eqn:Dirichlet:cts} also satisfies the estimate~\eqref{eqn:N:bound}.

Then there is a number $\mathfrak{b}^*\in (0,1]$ such that the conclusions of Theorem~\ref{thm:Dirichlet} are valid for $(s,1/p)$ in the region indicated on the left in Figure~\ref{fig:extrapolation}.

If instead there is a $q^*\in (1,\infty)$ such that the solution $u$ to the problem~\eqref{eqn:Dirichlet:cts}, with $A$ replaced by~$A^*$, satisfies the estimate~\eqref{eqn:N:bound} (with $q$ replaced by~$q^*$), then there is a $\mathfrak{b}\in (0,1]$ such that the conclusions of Theorem~\ref{thm:Dirichlet} are valid for $(s,1/p)$ in the region indicated on the right in Figure~\ref{fig:extrapolation}.

Furthermore, the numbers $\mathfrak{b}^*$ or $\mathfrak{b}$ satisfy the estimates
\begin{equation}\label{eqn:b:lower}
\mathfrak{b}\geq \alpha,
\qquad
\mathfrak{b}^*\geq \alpha^*,
\qquad
\mathfrak{b}\geq 1+\frac{\bdmn}{q^*}-\bdmn,
\qquad
\mathfrak{b}^*\geq 1+\frac{\bdmn}{q}-\bdmn
\end{equation}
where $\alpha$ and $\alpha^*$ are the boundary De Giorgi-Nash exponents for $L$ and~$L^*$.

(In the event that $q$ and $q^*$ both exist, the conclusions of Theorem~\ref{thm:Dirichlet} are valid for $(s,1/p)$ in the region indicated on the bottom in Figure~\ref{fig:extrapolation}.)
\end{thm}

Theorem~\ref{thm:Dirichlet} is essentially the $q=q^*=\infty$ case of Theorem~\ref{thm:Dirichlet:Lp}; note that the estimate~\eqref{eqn:N:bound} for $q=\infty$ reduces to the classical maximum principle.
Expressed algebraically, the pair $(s,1/p)$ satisfies the conditions of Theorem~\ref{thm:Dirichlet} if
\begin{align}
\label{eqn:thm:Dirichlet:Lq:infty}
p&=\infty&&&\text{and}&& 0&<s<\mathfrak{b},\\
\label{eqn:thm:Dirichlet:Lq:quasi}
p&\leq 1 &&&\text{and}&& \frac{\bdmn}{p}-\bdmn+1-\mathfrak{b}^*&<s<1,\\
\noalign{\noindent\text{or}}
\label{eqn:thm:Dirichlet:Lq:Banach}
1<p&<\infty,&& 0<s<1,&\text{and}&&
\frac{q-p}{qp-p}(1-\mathfrak{b}^*)
&<s
<\mathfrak{b}+\frac{q^*}{pq^*-p}(1-\mathfrak{b})
.\end{align}

As in Remark~\ref{rmk:convex}, the condition $\mathfrak{b}^*\geq 1+\bdmn/q-\bdmn$ implies that the regions in Figure~\ref{fig:extrapolation} are convex in a neighborhood of the point $(1-\mathfrak{b}^*,1)$, and thus are convex.

Our methods of proof do not require that the number $\bdmn$ in Theorem \ref{thm:Dirichlet} or~\ref{thm:Dirichlet:Lp} be an integer. However, in most of the cases  in which the nontangential estimate~\eqref{eqn:N:bound} is known to hold for a finite~$q$, the boundary $\partial\Omega$ is of integer dimension. (Some exceptions do exist; see in particular \cite{DavEM21}.) That is, in most of the known applications of Theorem~\ref{thm:Dirichlet:Lp}, $\bdmn$ is an integer.

\section{Further definitions}\label{sec:dfn:2}

In Section~\ref{sec:Lp} we will discuss some known cases in which the $L^q$-Dirichlet problem is well posed, and thus some cases in which we may apply Theorem~\ref{thm:Dirichlet:Lp}. It is known that the $L^q$-Dirichlet problem requires additional conditions on the coefficients and domains, beyond those listed in Theorem~\ref{thm:Dirichlet}; see \cite{ModM80,CafFK81}. That is, such results necessarily require more stringent conditions on the coefficients and domains. Thus, in this section we will define some of those more stringent conditions.

\subsection{Geometric conditions}

\begin{defn}[Uniformly rectifiable]\label{dfn:UR}Let $\bdmn$, $\dmn$ be integers with $1\leq \bdmn\leq \dmnMinusOne$, and let $F$ be a closed subset of $\R^\dmn$. We say that $F$ is uniformly rectifiable (or uniformly $\bdmn$-rectifiable) if $F$ is $\bdmn$-Ahlfors regular and has big pieces of Lipschitz images; that is, if there exist positive constants $\theta$ and $M$ such that if $x\in F$ and $0<r<\diam F$, then there is a Lipschitz function $\varphi:B_\bdmn(0,r)\to\R^\dmn$, where $B_\bdmn(0,r)\subset\R^\bdmn$ is the $\bdmn$-dimensional ball of radius~$r$, such that the Lipschitz constant of $\varphi$ is at most $M$ and such that
\begin{equation*}\sigma\bigl(F\cap B_\dmn(x,r)\cap \varphi(B_\bdmn(0,r))\bigr) \geq \theta r^\bdmn.\end{equation*}
\end{defn}

\begin{defn}[Harnack chains]\label{dfn:Harnack}
Let $\Omega\subset\R^\dmn$. We say that $\Omega$ satisfies the \emph{Harnack chain condition} if there are constants $M\geq 1$ and $m\in (0,1)$ such that if $x$, $y\in \Omega$, then there is an integer~$N$, depending only on the ratio $|x-y|/\delta$ for $\delta=\min(\dist(x,\partial\Omega),\dist(y,\partial\Omega))$, such that there is a chain of open balls $\{B_k\}_{k=1}^N$ with $x\in B_1$, $y\in B_N$, $B_k\cap B_{k+1}\neq\emptyset$ for all $1\leq k<N$, and $m\diam(B_k)\leq \dist(B_k,\partial\Omega)\leq M\diam(B_k)$.

\end{defn}

By \cite[Lemma~2.1]{DavFM21}, and as noted in Remark~\ref{rmk:higher:John:corkscrew}, if $\partial\Omega$ is $\bdmn$-Ahlfors regular for some $0<\bdmn<\dmnMinusOne$, then $\Omega$ automatically satisfies the Harnack chain condition. However, the higher codimension condition does not guarantee uniform rectifiability.

\begin{defn}[Chord-arc domain]\label{dfn:NTA}
Let $\Omega\subset\R^\dmn$. We call $\Omega$ a chord-arc domain if $\Omega$ satisfies the interior corkscrew condition (Definition~\ref{dfn:iCS}) and the Harnack chain condition, if $\partial\Omega$ is $\pdmnMinusOne$-Ahlfors regular, and if $\Omega$, in addition, satisfies the exterior corkscrew condition, that is, if the exterior domain $\R^\dmn\setminus\overline\Omega$ also satisfies the interior corkscrew condition.
\end{defn}

\begin{rmk} By connecting corkscrew points with Harnack chains, we can construct a $\lambda$-carrot path from any point $x\in\partial\Omega$ to any point $y\in \Omega$ with $|x-y|<N\dist(y,\partial\Omega)$, where $\lambda$ depends only on~$N$, the dimension~$\dmn$, and the corkscrew and Harnack chain constants. Thus, a chord-arc domain is a weak local John domain, and so satisfies the conditions of Theorem~\ref{thm:Dirichlet}.
\end{rmk}

\begin{rmk} A chord-arc domain necessarily has uniformly rectifiable boundary; this follows from \cite{DavJ90} or \cite{Sem90B} combined with \cite{HofMU14}. \end{rmk}

\begin{defn}[Chord-arc surfaces with small constant]\label{dfn:CASSC} Let $\Gamma\subset\R^{\dmn}$ be closed and unbounded. We say that $\Gamma$ is a chord-arc surface with constant~$\varepsilon$ if, for every $x\in \Gamma$ and $r>0$, there is a Lipschitz function $\varphi_x:\R^\dmnMinusOne\to\R$ with Lipschitz constant at most $\varepsilon$ and a rigid motion $T_x$ such that, if
\begin{equation*} G=\{T_x((y,\varphi(y))):y\in\R^\dmnMinusOne\}, \end{equation*}
then $G\cap B(x,r/2)\neq \emptyset$ and
\begin{equation*}\sigma(\Gamma\cap B(x,r)\setminus G)+\sigma(G\cap B(x,r)\setminus \Gamma)<\varepsilon r^\dmnMinusOne.\end{equation*}
\end{defn}

\begin{rmk}\label{rmk:CASSC}By \cite[Lemma~5.4]{DavLM23} and the following remarks, if $\Omega\subset\R^{\dmn}$ and $\partial\Omega$ is a chord-arc surface with sufficiently small constant~$\varepsilon$, then $\Omega$ and $\R^{\dmn}\setminus\overline\Omega$ are both chord-arc domains. \end{rmk}

\subsection{Conditions on coefficients}

In this section we will state one form of the Dahlberg-Kenig-Pipher (DKP) condition on coefficients~$A$. Early forms of this condition were proposed by Dahlberg and investigated by Kenig and Pipher in \cite{KenP01}, and many weaker forms have been studied.

In particular, for simplicity of exposition we will state the DKP condition only in the classical codimension~$1$ case (that is, $\bdmn=\dmnMinusOne$). However, generalizations of these notions to the higher codimensional case have been investigated (see \cite{DavFM19A,FenMZ21}) and are also a present active area of research.

\begin{defn}[The weak DKP condition]\label{dfn:DKP}
Let $\Omega\subset\R^{\dmn}$ be open and such that $\partial\Omega$ is $\pdmnMinusOne$-Ahlfors regular.
Let $A$ be a matrix-valued function defined on $\Omega$ that satisfies the ellipticity conditions~(\ref{eqn:elliptic:introduction}--\ref{eqn:elliptic:bounded:introduction}) for some $\Lambda>\lambda>0$ with $\bdmn=\dmn-1$.

We say that $A$ satisfies the weak DKP condition with constant $M$ in $\Omega$ if $A=B+C$, where
\begin{equation*}\sup_{x\in\partial\Omega}\sup_{r>0} \frac{1}{r^\dmnMinusOne} \int_{B(x,r)\cap\Omega} |\nabla B(y)|^2\dist(y,\partial\Omega) + |C(y)|^2\frac{1}{\dist(y,\partial\Omega)} \,dy =M<\infty.\end{equation*}
\end{defn}

\section{Corollaries of Theorem~\ref{thm:Dirichlet:Lp}}\label{sec:Lp}

The $L^q$-Dirichlet problem has been a subject of intense study over the past few decades, and so the bound~\eqref{eqn:N:bound} is known to be valid in many cases beyond Lipschitz domains. A complete survey of such results is beyond the scope of this note; in this section, we will describe a few such results.

\subsection{Higher codimensional results}

The most elegant (albeit not the most general) of the known results for the Dirichlet problem in the higher codimensional case was established in
\cite{DavM22} (and later by another method in \cite{Fen20p}) and is as follows.

Let $\bdmn$ be an integer with $1\leq\bdmn\leq \dmn-2$. Let $\Omega\subset\R^\dmn$ and suppose that $\partial\Omega$ is $\bdmn$-Ahlfors regular and uniformly rectifiable. (By Remark~\ref{rmk:higher:John:corkscrew}, because $\bdmn<\dmnMinusOne$, we have that $\Omega$ satisfies the conditions of Theorem~\ref{thm:Dirichlet}.)

If $\beta>0$, define the function $D_\beta:\Omega\to\R$ by
\begin{equation*}D_\beta(x)=\biggl(\int_{\partial\Omega} |x-y|^{-\bdmn-\beta}\,d\sigma(y) \biggr)^{-1/\beta}.\end{equation*}
By \cite[Lemma~5.1]{DavFM19A}, $D_\beta(x)\approx \dist(x,\partial\Omega)$.

Let $A_\beta(x)=D_\beta(x)^{\bdmn+1-\pdmn}I$, where $I$ is the identity matrix, and let $L_\beta=\Div A_\beta\nabla$ (as defined by Definition~\ref{dfn:weak}). Observe that $A$ satisfies the ellipticity conditions (\ref{eqn:elliptic:introduction}--\ref{eqn:elliptic:bounded:introduction}), and also is a symmetric matrix and so $L_\beta=L_\beta^*$.

By \cite[Theorem~1.14]{Fen20p} or \cite[Theorem~1.1]{DavM22}, the harmonic measure associated to $L$ is in $A_\infty(\sigma)$; by \cite[Theorem~4.1]{MayZ19}, this implies that the $L^q$-Dirichlet problem is well posed for some $q\in (1,\infty)$.

Thus, for a uniformly rectifiable boundary, and for the operator $L_\beta$, $\beta>0$, we have that the conditions of Theorem~\ref{thm:Dirichlet:Lp} are satisfied, and so the solution to the Poisson-Dirichlet problem~\eqref{eqn:Dirichlet:besov} satisfies the estimate~\eqref{eqn:Besov:estimate} for all $(s,1/p)$ in the hexagon at the bottom of Figure~\ref{fig:extrapolation}.

Especially in the higher codimensional setting, the question of coefficients $A$ for which the $L^q$-Dirichlet problem is well posed is an active area of research and we expect further results of this type to be proven soon. Some (somewhat complicated) results of this type are known; we refer the interested reader to \cite{DavFM19A,FenMZ21} for the details.

\subsection{The Laplace operator and the DKP condition}

We now turn to known results for the codimension~$1$ case, that is, the case $\bdmn=\dmn-1$.

As in the previous section, the most elegant known results are those for a particular operator; as is usual in the codimension one case, the favored particular operator is the Laplacian.
By \cite{AzzHMMT20}, if $\Omega\subset\R^\dmn$ is a weak local John domain and satisfies the interior corkscrew condition, and if $\partial\Omega$ is $\pdmnMinusOne$-Ahlfors regular and uniformly rectifiable, then there is a $q\in (1,\infty)$ such that the $L^q$-Dirichlet problem is well posed for $L=-\Delta$ (that is, for coefficients $A=I$, where $I$ is the identity matrix).

There has been considerable investigation into the weakest possible sufficient conditions (and the strongest possible necessary conditions!) for well posedness of the $L^q$-Dirichlet problem. Many forms of the DKP condition have been investigated under many conditions on domains. We will now mention one known well posedness result for operators satisfying a quite weak version of the DKP condition (defined above); this is not the only such known result and in fact some results are known in more general domains.

Specifically, if $\Omega\subset\R^\dmn$ is a chord-arc domain, and if $A$ satisfies the weak DKP condition of Definition~\ref{dfn:DKP} in~$\Omega$ for some finite~$M$, then there is a $q\in (1,\infty)$ such that the $L^q$-Dirichlet problem for $-\Div A\nabla$ is well posed.

This follows from known results as follows. By \cite[Theorem~1.21]{FenLM24}, if $\Omega\subset\R^\dmn$ is open, if $\partial\Omega$ is $\pdmnMinusOne$-Ahlfors regular and uniformly rectifiable, if $A$ satisfies the weak DKP condition, and if $-\Div A\nabla u=0$ and $u$ is bounded in~$\Omega$, then
\begin{equation}\label{eqn:FenLM24}
\int_{B(x,r)\cap\Omega} |\nabla u(y)|^2\dist(y,\partial\Omega)\,dy\leq C \|u\|_{L^\infty}^2 \sigma(B(x,r)\cap\partial\Omega)
\end{equation}
for all $x\in\partial\Omega$ and all $r>0$.

By \cite[Theorem~1.1]{CavHMT20}, if $\Omega$ is a one-sided chord-arc domain, then the Carleson condition~\eqref{eqn:FenLM24} implies that the harmonic measure for $-\Div A\nabla$ lies in $A_\infty(\partial\Omega)$.

Let $\omega^y$ be the harmonic measure for $-\Div A\nabla$ with pole at~$y$, and let $k^y$ be the Radon-Nikodym derivative of $\omega^y$ with respect to the surface measure.
As observed just following \cite[Definition~2.13]{CavHMT20}, their $A_\infty$ condition implies that there exists a $C<\infty$ and a $p\in (1,\infty)$ such that, if $x\in \partial\Omega$, $0<r<\diam\Omega$, and $y=A_r(x)$ is the corkscrew point to $x$ at scale~$r$, then
\[\biggl(\fint_{\Delta(x,r)} (k^{y})^p\,d\sigma\biggr)^{1/p}
\leq C\fint_{\Delta(x,r)} k^{y}\,d\sigma.\]

It is well known (see, for example, \cite[Proposition~4.5]{HofL18}) that this reverse H\"older estimate on the harmonic measure yields well posedness of the $L^q$-Dirichlet problem for $q=p/(p-1)$, that is, for $q$ the H\"older conjugate of~$p$.

\subsection{The small DKP condition}
We conclude this discussion by mentioning a ``small constant'' result. That is, notice that the results of the previous two sections all imply existence of a $q\in (1,\infty)$ such that the $L^q$-Dirichlet problem is well posed, without ever specifying precise values of~$q$. By requiring more of our coefficients and domains, it is possible to choose a value of~$q$ and find situations in which the $L^q$-Dirichlet problem is well posed.

Specifically, let $q\in (1,\infty)$, let $\dmn\geq 3$ be an integer, and let $\Lambda>\lambda>0$. Then there is a $\delta>0$ (depending on $q$, $\bdmn$, $\Lambda$, and~$\lambda$) such that if
\begin{enumerate}
\item $\Omega\subset\R^{\dmn}$,
\item $\partial\Omega$ is a chord-arc surface with constant at most~$\delta$ (as defined in Definition~\ref{dfn:CASSC}),
\item $A$ satisfies the ellipticity conditions~(\ref{eqn:elliptic:introduction}--\ref{eqn:elliptic:bounded:introduction}) (with $\bdmn=\dmn-1$), and
\item $A$ satisfies the weak Dahlberg-Kenig-Pipher condition in $\Omega$ with constant at most~$\delta$ (as defined in Definition~\ref{dfn:DKP}),
\end{enumerate}
then the $L^q$-Dirichlet problem is well posed.

This result follows from known results in the literature by a line of argument similar to the previous section. First, by the recent results of \cite{DavLM23}, the given conditions imply a quantitative local $A_\infty$-type condition on the harmonic measure.\footnote{The weak $DKP$ condition required by \cite{DavLM23} is not quite that required by \cite{FenLM24}. To lighten the exposition, we have chosen to state only the stronger of these two weak $DKP$ conditions. It is straightforward to establish that the weak $DKP$ condition of \cite{DavLM23} implies the weak $DKP$ condition of \cite{FenLM24} in the setting of \cite{DavLM23}, that is, in domains $\Omega$ such that the Whitney boxes of \cite[formula~(1.2)]{DavLM23} have volume comparable to~$r^\dmn$.} Second, one can show that this $A_\infty$-type condition implies a local reverse H\"older condition on the density (that is, the Radon-Nikodym derivative with respect to the surface measure) of the harmonic measure. Finally, as before, by \cite[Proposition~4.5]{HofL18}, this reverse H\"older estimate yields well posedness of the $L^q$-Dirichlet problem.

The precise form of the $A_\infty$-type condition varies from paper to paper and some work must be done to relate particular forms to the required reverse H\"older condition. Moreover, in the small constant case, we must at every step ensure that our results work for the chosen value of~$q$, and not, as in the previous section, merely for sufficiently large finite~$q$. This means that the argument for this result is somewhat more involved than in the previous section; we will give the details in the appendix to \cite{BarMP25pA}.

\medskip

We emphasize that all of the results mentioned in Section~\ref{sec:Lp} (and many other recent advancements in this spirit) extrapolate automatically to well-posedness within the subregions in Figure~\ref{fig:extrapolation} by virtue of Theorem~\ref{thm:Dirichlet:Lp}. In particular, in the small constant case, the region includes arbitrarily large compact subsets of the square $(0,1)\times (1, \infty)$. We emphasize that this is not a result of the more traditional interpolation between the Dirichlet problem and the regularity problem, as in most of these cases the regularity problem has not been treated.

\def\cprime{'}
\bibliographystyle{amsalpha}
\bibliography{bibli}
\end{document}